\documentclass[ejs,preprint]{imsart}

\usepackage{amsthm,amsmath}
\RequirePackage[numbers]{natbib}
\RequirePackage[colorlinks,citecolor=blue,urlcolor=blue]{hyperref}

\doi{}
\pubyear{0000}
\volume{0}
\firstpage{0}
\lastpage{0}

\startlocaldefs

 \usepackage[utf8]{inputenc}
 \usepackage{amsmath}
 \usepackage{amsfonts}
 \usepackage{amssymb}
\usepackage{amsthm}

 \usepackage{dsfont}
 \usepackage{bm}

 \usepackage{graphicx}
 \usepackage{enumerate}
 \usepackage{paralist}

 \newtheorem{theorem}{Theorem}[section]
\newtheorem{lemma}[theorem]{Lemma}
\newtheorem{proposition}[theorem]{Proposition}
 \newtheorem{cor}[theorem]{Corollary}

\newenvironment{proofof}{\noindent\sc Proof of}{
    \hspace*{\fill} $\Box$ \vspace{2ex} }

\newcommand{\cond}{\textbf{C}\!\!}

\theoremstyle{remark}
\newtheorem{remark}[theorem]{Remark}
\numberwithin{equation}{section}

\renewcommand{\P}{\mathbb{P}}
\newcommand{\E}{\operatorname{\mathbb{E}}}

\newcommand{\cov}{\operatorname{\mathbf{Cov}}}
\newcommand{\e}{\operatorname{e}}

\DeclareMathOperator{\tr}{tr}
\DeclareMathOperator{\Span}{span}
\DeclareMathOperator{\supp}{supp}

\usepackage{xspace}

\newcommand\cf{\emph{cf.}\xspace }

\newcommand{\inprod}[1]{\langle #1 \rangle}

\newcommand{\sphere}{\mathbb{S}^{d-1}}

\newcommand{\proj}{\mathbf{\Pi}}
\newcommand{\ang}{\Theta}
\newcommand{\angf}{\theta}
\newcommand{\HSnorm}[1]{{\left\| #1\right\|_{HS}^2}}
\newcommand{\hsnorm}[1]{{\left\| #1\right\|_{HS}}}
\newcommand{\opnorm}[1]{{\left\| #1\right\|_{op}}}
\newcommand{\hsprod}[2]{{\left\langle #1,#2\right\rangle_{HS}}}

\newcommand{\R}{\mathbb{R}}
\newcommand{\C}{\mathbb{C}}
\newcommand{\N}{\mathbb{N}}

\newcommand{\Z}{\mathbb{Z}}
\newcommand{\M}{\mathbb{M}}
\newcommand{\eps}{\varepsilon}
\newcommand{\VV}{\mathcal{V}}
\newcommand{\MM}{\mathcal{M}}
\newcommand{\FF}{\mathcal{F}}
\newcommand{\HH}{\mathcal{H}}
\newcommand{\WW}{\mathcal{W}}
\newcommand{\II}{\mathcal{I}}
\newcommand{\RR}{\mathcal{R}}
\newcommand{\LL}{\mathcal{L}}
\newcommand{\BB}{\mathcal{B}}
\newcommand{\Ub}{\mathbf{U}}
\newcommand{\Wb}{\mathbf{W}}
\newcommand{\Ind}[1]{{\boldsymbol{1}_{\{#1\}}}}
\newcommand{\ind}[1]{{\boldsymbol{1}_{#1}}}

\usepackage[normalem]{ulem}

\endlocaldefs

\allowdisplaybreaks

\begin{document}

\begin{frontmatter}

\title{Asymptotic Behavior of Principal Component Projections for Multivariate Extremes} 
\runtitle{Asymptotics for PCA for Extremes, \today}


\author{\fnms{Holger} \snm{Drees}\ead[label=e1]{holger.drees@uni-hamburg.de}}
\address{  University of Hamburg, Department of Mathematics, Germany\\
  \printead{e1}}

\runauthor{H. Drees}

\begin{abstract}

The extremal dependence structure of a regularly varying $d$-dimensional random vector can be described by its angular measure. The standard nonparametric estimator of this measure is the empirical measure of the observed angles of the $k$ random vectors with largest norm, for a suitably chosen number $k$. Due to the curse of dimensionality, for moderate or large $d$, this estimator is often inaccurate. If the angular measure is concentrated on a vicinity of a lower dimensional subspace, then first projecting the data on a lower dimensional subspace obtained by a principal component analysis of the angles of extreme observations can substantially improve the performance of the estimator.

We derive the asymptotic behavior of such PCA projections and the resulting excess risk. In particular, it is shown that, under mild conditions, the excess risk (as a function of $k$) decreases much faster than it was suggested by empirical risk bounds obtained in \cite{DS21}. Moreover, functional limit theorems for local empirical processes of the (empirical) reconstruction error of projections uniformly over neighborhoods of the true optimal projection  are established. Based on these asymptotic results, we propose a data-driven method to select the dimension of the projection space. Finally, the finite sample performance of resulting estimators is examined in a simulation study.

\end{abstract}

\begin{keyword}[class=MSC2020]
\kwd[Primary ]{62G32}
\kwd[; secondary ]{62H25, 62G20}
\end{keyword}
\begin{keyword}[class=KWD]
\kwd{   dimension reduction,  multivariate extreme value analysis,  multivariate regular variation, Principal Component Analysis, projection matrix, local empirical processes}
\end{keyword}



\end{frontmatter}

\section{Introduction}
\label{sec:intro}

Suppose that we want to analyze the extreme value behavior of a random vector $X\in\R^d$. For example, the components of $X$ may describe the returns of different financial assets and we are interested in the Value at Risk (i.e.\ a high quantile) of an arbitrary portfolio built from these assets. In environmetrics, the components of $X$ may model flood heights at different points along a coast protected by a dike, where a crevasse at any point will lead to a flooding of the whole area behind the dike. If one is interested in probabilities of events which have rarely (or never) happened yet, then it is necessary to work in a suitable extreme value framework.

\subsubsection*{The setting}
A standard assumption is that $X$ is regularly varying: there exists a non-degenerate measure $\mu$ on $(\R^d\setminus\{0\},\BB(\R^d\setminus\{0\}))$ such that for all $\mu$-continuity sets $B\in\BB(\R^d\setminus\{0\})$ bounded away from 0 one has
$$ \lim_{t\to\infty} \frac{\P\{X/t \in B\}}{\P\{\|X\|>t\}} = \mu(B)<\infty. $$
An equivalent formulation is that the conditional distribution of $X$ given that $\|X\|>t$ converges weakly to $\mu$ restricted to $\{x\in\R^d | \|x\|>1\}$.
Here $\|\cdot\|$ could be any norm on $\R^d$ (changing the norm merely alters $\mu$ by a constant factor), but in what follows we will always consider the Euclidean norm. Then automatically $\mu$ is homogeneous of order $-\alpha$, for some $\alpha>0$, the so-called index of regular variation. Moreover, the survival function of $\|X\|$ is regularly varying, too, i.e.\ $\lim_{t\to\infty}\P\{\|X\|>tx\}/\P\{\|X\|>t\}=x^{-\alpha}$ for all $x>0$. (For that reason, in environmental applications one must usually first standardize the marginal distributions to make regular variation a realistic assumption.)
The homogeneity of $\mu$ implies that it allows the following spectral decomposition. For the so-called angular (or spectral) probability measure $H$ on $(\sphere,\BB(\sphere))$, defined by
$$ H(A) := \mu\Big\{x\in\R^d\,\Big|\, \|x\|>1, \frac{x}{\|x\|}\in A\Big\}, $$
and the polar transformation
$$ T:\R^d\setminus\{0\}\to (0,\infty)\times\sphere, \; T(x)=\Big(\|x\|,\frac{x}{\|x\|}\Big) $$
the push-forward measure $\mu^T$ equals the product of a the measure on $(0,\infty)$ with density $x\mapsto\alpha x^{-(\alpha+1)}$  and $H$. So while, by and large, the marginal extreme value behavior is determined by the single parameter $\alpha$, the extreme value dependence of $X$ is governed by the angular measure, which can be any probability measure on $\sphere$.

\subsubsection*{Reduction of complexity}
Therefore, multivariate extreme value theory is inherently non-parametric and thus prone to the curse of dimensionality. To render standard methods of multivariate extreme value theory practically feasible in settings of high (and often even moderate) dimension, some reduction of the model complexity is needed. Several different approaches to this goal have been suggested. For example, the extreme value dependence between the components of the vector can be described by a graph structure (with the components as nodes). If this graph has few edges, then the dependence structure is much simplified; see e.g.\ \cite{EH20}. Alternatively, one may try to split the vector in subvectors of relatively low dimension which are asymptotically independent in extreme regions; cf.\ \cite{CSS19}. See \cite{EI21} for a nice review of methods for complexity reduction techniques in extreme value theory.

Here we discuss a version of the classical principal component analysis (PCA) suggested by \cite{DS21}, which is tailored to our extreme value setting. Throughout, we assume that $n$ independent and identically distributed (iid) random vectors $X_i$, $1\le i\le n$, are observed that are regularly varying with limit measure $\mu$ and index $\alpha$. As the components of $X$ have infinite second moment if $\alpha<2$, one must first standardize the vector in a suitable way before one can apply PCA. For that reason, we define
\begin{equation} \label{eq:def_ang}
  \ang_i := \omega(X_i) X_i, \quad 1\le i\le n,
\end{equation}
where $\omega:\R^d\to (0,\infty)$ is some measurable scaling function that ensures that $\ang_i$ is bounded. By far the most common choice will be $\omega(x)=1/\|x\|$ for $x\ne 0$, so that the conditional distribution of $\ang_i\in\sphere$ given that $\|X_i\|>t$ converges to the angular measure $H$. Since we are interested in the extreme value behavior described by the angular measure, we discard all observations with small or moderate norm and apply PCA to the remaining re-scaled vectors. Note that such a simple re-scaling of the observations ensures that the findings from PCA applied to those vectors $\ang_i$ for which $\|X_i\|$ is large  have an immediate  interpretation for the original observations. In contrast, \cite{CT19} (see also \cite{JCW20}) proposed first to  transform each coordinate of the vector non-linearly and then use PCA to analyze some estimator of the so-called tail pairwise dependence matrix. Since the pre-images of linear subspaces under the (by and large arbitrarily chosen) non-linear transformation are usually non-linear manifolds, the interpretation of the PCA approximation and the reconstruction error for the original data is difficult. In \cite{DS21} finite sample bounds on the excess risk are are derived for the present setting. These results were generalized to PCA for extremes of Hilbert space values random variables in \cite{CHS24}.

Another related article is \cite{Chautru15}, in which the observations are re-scaled to the unit sphere, but then successively sub-spheres of decreasing dimension are determined  which have the smallest reconstruction error w.r.t.\ the spherical distance. Since these sub-spheres in general are not intersections of linear subspaces with the sphere, the resulting decomposition of the data is fundamentally different from the linear projections analyzed here.
\cite{AMDS22} applies a kernel version of PCA to extremes, which are first mapped by a feature map into an infinite dimensional Hilbert space. The goal, though, is very different, namely not to estimate the angular measure accurately, but to detect clusters in which the angular measure is approximately concentrated, while any mass of the angular measure in other regions (which is considered noise) should be pushed to these clusters.  Since these alternatives are substantially more involved (in case of \cite{Chautru15} also numerically), no mathematical results on their performance in terms of the reconstruction error are available. Other approaches to clustering for the angular measure were suggested by \cite{JW20} and \cite{AMDS24}.

\subsubsection*{Outline of PCA for extremes}
To give a more precise description of the PCA analyzed here, denote by $R_i:=\|X_i\|$ the norm of the $i$th observation and, for some $k\in\{1,\ldots,n\}$, the $k+1$ largest order statistic of $R_i$, $1\le i\le n$, by $\hat t_{n,k}$. For a fixed $p\in\{1,\ldots,d-1\}$, we determine a $p$-dimensional linear subspace $V$ with projection matrix $\proj_V$ such that
\begin{equation} \label{def:hatMnk}
 \hat M_{n,k}(\proj_V) := \frac 1k \sum_{i=1}^n \ang_i^\top \proj_V \ang_i \Ind{R_i>\hat t_{n,k}}
\end{equation}
is maximized. This requirement is equivalent to the condition that the (empirical) reconstruction error $\|\ang_i-\proj_V\ang_i\|^2$ summed over the $k$ observations with largest norm is minimal. Denote this subspace by $\hat V_{n,k}$ and its theoretical counterpart that maximizes
\begin{equation} \label{def:Mnk}
   M_{n,k}(\proj_V) := \E(\ang^\top \proj_V \ang\mid R>t_{n,k})
\end{equation}
by $V_{n,k}^*$, with $t_{n,k}$ denoting the $(1-k/n)$-quantile of $R:=\|X\|$. In \cite{DS21}  finite sample uniform bounds on $|\hat M_{n,k}(\proj_V)-M_{n,k}(\proj_V)|$ have been established. In particular, they imply (for fixed dimensions $p$ and $d$) an upper bound on the so-called excess risk $M_{n,k}(\proj_{V_{n,k}^*})-M_{n,k}(\proj_{\hat V_{n,k}})$ of the order $k^{-1/2}$. It turns out, however, that this bound is  grossly misleading, because in fact the excess risk decreases with the rate $k^{-1}$ if there is a unique maximizer in the limit model. The effect (which is known from classical PCA) is easily explained by the different behavior of the functional $M_{n,k}$ near its maximum and in a neighborhood of some other projection matrix $\proj_V$.

The goals of the present article is three-fold. First, we examine the asymptotic behavior of the PCA projection matrix and the resulting excess risk. In particular,  we derive limit distributions for these quantities in the Corollaries \ref{cor:limitPCAproj} and \ref{cor:excessriskunique} if the maximizer in the limit model is unique. If this is not the case, then one should not expect that the PCA projections converge. However, one can still prove rates of convergence on the excess risk which will often be faster than $k^{-1/2}$. In Section 5, we analyze both the empirical risk and its theoretical counterpart for projections in a neighborhood of the optimal projections. This analysis helps to  better understand the asymptotic behavior of these risk functions in the relevant regions and yields a different interpretation of the respective maximizers.

These results are then used to suggest a new method for selecting the dimension $p$ of the subspaces on which we project. In this context, it is important to understand that, unlike in the classical setting, there may be a true dimension $p<d$ which is non-trivial to estimate. If in the classical setting the ground distribution of the observations is concentrated on a $p$-dimensional subspace with $p<d$, then this would be obvious from the data. In the present situation, though, if the limit measure $\mu$ is concentrated on such a subspace, in general even the observations with large norm will still span the complete $\R^d$, because they are drawn from an intermediate distribution which will usually not be concentrated on the same set as the limit measure.

\subsubsection*{Notation}

All random variables are defined on some probability space $(\mathcal{X},\mathcal{A},\P)$. The expectation with respect to $\P$ is denoted by $\E$, and the law of some random element $Y$ by $\LL(Y)$.
For an $\R^d$-valued random vector $X$ and a scale function $\omega$ as above, let $R=\|X\|$ be its Euclidean norm, $\ang=\omega(X)X$ and let $\ang_t$ be a random vector with the same distribution as $\ang$ conditionally on $R>t$.
$\E_t$ denotes the conditional expectation given $R>t$, and $\E_\infty$ the expectation in the limit model, that is w.r.t.\ $\mu$ restricted to the complement of the unit ball. The $(1-k/n)$-quantile of $R$ is denoted by $t_{n,k}$, and $\hat t_{n,k}$ denotes its empirical counterpart, i.e.\ the $k+1$ largest order statistics of the random variables $R_i=\|X_i\|$, $1\le i\le n$.

The transposition of a vector or a matrix is denoted by the superscript $^\top$. For any  subset $D$ of $\R^d$ (or family in $\R^d$), $\Span(D)$ denotes the linear subspace spanned by the elements of $D$.
For an arbitrary linear subspace $V$ of $\R^d$, $V^\perp$ denotes its orthogonal complement and $\proj_V$ the projection (matrix) onto $V$. The set of all $p$-dimensional subspaces of $\R^d$ is denoted by $\VV_p$. The operator norm of a matrix $A\in\R^{d\times d}$ is denoted by $\|A\|_{op}:=\sup_{x\in \sphere} \|Ax\|$, with $\sphere$ denoting the unit sphere in $\R^d$. The Hilbert-Schmidt product of two matrices $A,B\in \R^{d\times d}$ is given by $\hsprod{A}{B} = \tr(AB^\top)$, with $\tr$ denoting the trace functional. The corresponding matrix norm is denoted by $\hsnorm{\cdot}$. Since all matrix norms are equivalent (if the dimension is fixed), rates of convergence may be given in any of them, but in different contexts it is more natural to work either with the operator or the Hilbert-Schmidt norm.


\section{Convergence of empirical mixed moments}
\label{sect:convempmom}

It is well known that any maximizer of $\hat M_{n,k}$ defined by \eqref{def:hatMnk} can be derived from the eigendecomposition of the  matrix of empirical mixed moments of the vectors $\ang_i$ with norms exceeding the random threshold $\hat t_{n,k}$, i.e.
$$ \hat\Sigma_{n,k} := \frac{1}{k} \sum_{i=1}^n \ang_i\ang_i^\top \Ind{R_i>\hat t_{n,k}}. $$
Its theoretical counterpart  is given by
$$ \Sigma_{n,k} := \E_{t_{n,k}}(\ang\ang^\top). $$
 Note that, by multivariate regular variation, $\Sigma_{n,k}$ converges to
$$ \Sigma_\infty := \E_\infty(\ang\ang^\top)
$$
w.r.t.\ the operator norm or the Hilbert-Schmidt norm. Since the PCA projection matrix can be written as a smooth function of $\hat\Sigma_{n,k}$, if it is unique, we first examine the asymptotic behavior of this matrix.

In what follows, we will often make the following assumption:\\[1ex]
{\bf (B)}\; \parbox[t]{11cm}{There exists $\eps>0$ such that for all $s\in[1-\eps,1+\eps]$
$$ \opnorm{\E_{st_{n,k}}(\ang\ang^\top)-\Sigma_{n,k}} = o(k^{-1/2}).$$}

This assumption may be considered a bias condition for the matrix of mixed second moments of  $\ang$ conditionally on the radius being large.  By regular variation, it is always fulfilled if $k$ does not grow too fast.
 Note, however, that we do not compare the matrix for a high threshold with  the limiting matrix $\Sigma_\infty$, but the matrices corresponding to exceedances over two high thresholds which differ only by a factor close to 1. Hence, one might expect that the conditional mixed moments are similar, which shows that this condition is very weak.

\begin{theorem} \label{theo:empmixedmatconv}
  If condition (B) is met, then for $\Delta_n:= \hat\Sigma_{n,k}-\Sigma_{n,k}$ one has
  $$ k^{1/2} \Delta_n \to \Ub \qquad \text{weakly,}$$
  where $\Ub=(U_{ij})_{1\le i,j\le d}$ is a symmetric centered Gaussian matrix with \\ $Cov(U_{ij},U_{\ell m}) =$  $Cov_\infty(\theta_i\theta_j,\theta_\ell\theta_m)$ for all $1\le i,j,\ell,m\le d$.
\end{theorem}
\begin{proof}
  Conditionally on $\hat t_{n,k}/t_{n,k}=s$, the empirical matrix $\hat\Sigma_{n,k}$ is distributed as $k^{-1}\sum_{i=1}^k \ang^{(n,s)}_i (\ang^{(n,s)}_i)^\top$, where $\ang^{(n,s)}_i$ are iid random vectors with distribution $\LL(\ang\mid R>st_{n,k})$. Denote a generic copy of these vectors by $\ang^{n,s}=(\theta_j^{(n,s)})_{1\le j\le d}$.

  Fix any $s\in[1-\eps,1+\eps]$.
  The regular variation of $X$ implies $Cov(\theta^{(n,s)}_i\theta^{(n,s)}_j,$ $\theta^{(n,s)}_\ell\theta^{(n,s)}_m)\to Cov_\infty(\theta_i\theta_j,\theta_\ell\theta_m)$.
  Since the random variables $\theta_i$ are bounded, a combination of the Cram\'{e}r-Wold device with the central limit theorem by Lindeberg and Feller shows that
  $$ k^{-1/2} \sum_{i=1}^k \Big(\ang_i^{(n,s)}(\ang_i^{(n,s)})^\top-\E_{st_{n,k}}(\ang\ang^\top)\Big) \to \Ub \quad \text{weakly.}
  $$
  Thus, Condition (B) yields
  $$ k^{-1/2} \sum_{i=1}^k \big(\ang_i^{(n,s)}(\ang_i^{(n,s)})^\top-\Sigma_{n,k}\big) \to \Ub \quad \text{weakly}
  $$
  for all $s\in[1-\eps,1+\eps]$. Because $\hat t_{n,k}/t_{n,k}\to 1$ in probability, by dominated convergence we may conclude
  \begin{align*}
  \Big| & \P\big\{k^{1/2}(\hat\Sigma_{n,k}-\Sigma_{n,k})\in B\big\}-\P\{\Ub\in B\}\Big| \\
  & \le \int_{1-\eps}^{1+\eps} \bigg| \P\Big\{  k^{-1/2} \sum_{i=1}^k \big(\ang_i^{(n,s)}(\ang_i^{(n,s)})^\top-\Sigma_{n,k}\big)
 \in B\Big\} -\P\{\Ub\in B\} \bigg|\\
  & \hspace{2.5cm} \,  \LL\Big(\frac{\hat t_{n,k}}{t_{n,k}}\Big)(ds)
  + \P\big\{\hat t_{n,k}/t_{n,k}\not\in[1-\eps,1+\eps]\big\}\\
 & \to 0,
  \end{align*}
  for all $\LL(\Ub)$-continuity sets $B\in\BB(\R^{d\times d})$, which proves the assertion.
\end{proof}

\begin{remark}
  \begin{enumerate}
     \item The proof shows that it would be sufficient to consider the supremum over all $s\in [1-c_{n,k},1+c_{n,k}]$ for any sequence $c_{n,k}$ such that $|\hat t_{n,k}/t_{n,k}-1|=o_P(c_{n,k})$ in Condition (B).
     \item One can also derive finite sample bounds on $\hat\Sigma_{n,k}-\Sigma_{n,k}$ if a uniform bound on the difference between $\E_{st_{n,k}}(\ang\ang^\top)$ and $\Sigma_{n,k}$ is known. Indeed, Theorem 9 of \cite{KoltLoun17} implies that there exists a constant $C$ such that for all $z>1$ one has
          \begin{align*}
           &\opnorm{\hat\Sigma_{n,k}-\Sigma_{n,k}} \le C \lambda_1^{(n)}\max\bigg(\sqrt{\frac{\max(t,1/\lambda_1^{(n)})}k},\frac{\max(t,1/\lambda_1^{(n)})}k\bigg)\\
           & \hspace{4cm} +\sup_{s\in[1-\eps,1+\eps]}
         \opnorm{\E_{st_{n,k}}(\ang\ang^\top)-\Sigma_{n,k}}
         \end{align*}
         with probability at least $1-\e^{-t}-\P\{|\hat t_{n,k}/t_{n,k}-1|>\eps\}$. Here $\lambda^{(n)}_1$ denotes the largest eigenvalue of $\Sigma_{n,k}$, which is at most 1 if in \eqref{eq:def_ang} we standardize with the norm. Using Bernstein's inequality, one can easily see that $\P\{|\hat t_{n,k}/t_{n,k}-1|>\eps\}$ decreases exponentially fast with $k$. (For an explicit bound on this term, in addition bounds on the deviation between the conditional distribution of the radius given that it exceeds a high threshold and its Pareto limit is needed.)
     \end{enumerate}
\end{remark}

\section{Convergence of projections onto eigenspaces}

By multivariate regular variation $\Sigma_{n,k}\to \Sigma_\infty$, and hence by Theorem \ref{theo:empmixedmatconv}
$\hat\Sigma_{n,k}\to \Sigma_\infty$ in probability. It is well known that this convergence implies the convergence of projections onto suitably matched eigenspaces. Indeed, we may even conclude that the difference between the projections multiplied with $\sqrt{k}$ converges weakly.

For a precise formulation of the result we need some additional notation. Denote by $\lambda_1^{(n)}\ge \lambda_2^{(n)}\ge\cdots\ge \lambda_d^{(n)}$, $\hat\lambda_1^{(n)}\ge \hat\lambda_2^{(n)}\ge\cdots\ge \hat\lambda_d^{(n)}$ and $\lambda_1^{(\infty)}\ge \lambda_2^{(\infty)}\ge\cdots\ge \lambda_d^{(\infty)}$ the eigenvalues of $\Sigma_{n,k}$, of $\hat\Sigma_{n,k}$,  and of $\Sigma_\infty$, respectively, in decreasing order. There exists an increasing sequence of numbers $N_j\in\{0,\ldots d\}$, $0\le j\le d_\infty$, such that $N_0=0$, $N_{d_\infty}=d$, $\lambda^{(\infty)}_m=\lambda^{(\infty)}_{N_j}=:\mu_j$ for all $N_{j-1}<m\le N_j$, $1\le j\le d_\infty$, and $\lambda_{N_j}^{(\infty)}>\lambda_{N_j+1}^{(\infty)}$ for all $1\le j<d_\infty$. Hence $\mu_j$, $1\le j\le d_\infty$, are the distinct eigenvalues of $\Sigma_\infty$ and $N_j$ is the maximum index of an eigenvalue equal to $\mu_j$.

Moreover, let $W_j^{(\infty)}$ be the eigenspace corresponding to $\mu_j$. Choose an orthonormal system of eigenvectors $v_m^{(n)}$, $1\le m\le d$, of $\Sigma_{n,k}$ corresponding to the eigenvalues in decreasing order and likewise an orthonormal system of eigenvectors $\hat v_m^{(n)}$, $1\le m\le d$, of $\hat\Sigma_{n,k}$. Finally, let $W_j^{(n)}=\Span\big((v_m^{(n)})_{N_{j-i}<m\le N_j}\big)$ and $\hat W_j^{(n)}=\Span\big((\hat v_m^{(n)})_{N_{j-1}<m\le N_j}\big)$.
\begin{theorem} \label{theo:eigenproj}
  If Condition (B) holds, then
  \begin{align}
    k^{1/2} & \big( \proj_{\hat W_j^{(n)}}- \proj_{W_j^{(n)}}\big)_{1\le j\le d_\infty} \nonumber \\
     &\to \bigg(\sum_{\substack{\ell=1\\ \ell\ne j}}^{d_\infty} \frac 1{\mu_\ell-\mu_j} \big(\proj_{W_\ell^{(\infty)}} \Ub \proj_{W_j^{(\infty)}} + \proj_{W_j^{(\infty)}} \Ub \proj_{W_\ell^{(\infty)}}\big)\bigg)_{1\le j\le d_\infty}
  \end{align}
  weakly.
\end{theorem}
\begin{proof}
  We employ the ideas of the proof of Theorem 5.1.4 in \cite{HE15}. Denote by $\RR_n(z) := (\Sigma_{n,k}-zI)^{-1}$ and $\hat\RR_n(z) := (\hat\Sigma_{n,k}-zI)^{-1}$, $z\in\C$, the resolvents of $\Sigma_{n,k}$ and $\hat\Sigma_{n,k}$, respectively. Fix some $j\in\{1,\ldots,d_\infty\}$ and let $\Gamma$ be a path along the circle in $\C$ with center $\mu_j$ and radius $\eta:= \min_{\ell\in\{1,\ldots,d_\infty\}\setminus\{j\}} |\mu_j-\mu_\ell| /2$. For sufficiently large $n$, one has $\max_{1\le m\le d} |\lambda_m^{(n)}-\lambda_m^{(\infty)}|<\eta/4$, so that $\Gamma$ encircles the eigenvalues $\lambda_m^{(n)}$ for $N_{j-1}<m\le N_j$, but no other eigenvalue of $\Sigma_{n,k}$. Likewise, if $\opnorm{\Delta_n}<\eta/2$, then $\max_{1\le m\le d} |\hat\lambda_m^{(n)}-\lambda_m^{(\infty)}|<\eta/4+\opnorm{\Delta_n}<3\eta/4$, and $\Gamma$ encircles the eigenvalues $\hat\lambda_m^{(n)}$ for $N_{j-1}<m\le N_j$, but no other eigenvalue of $\hat\Sigma_{n,k}$.

  According to Theorem 5.1.3 of \cite{HE15}, one has
  \begin{equation} \label{eq:projdiff1}
    \proj_{\hat W_j^{(n)}}- \proj_{W_j^{(n)}} = -\frac 1{2\pi i} \oint_\Gamma \hat\RR_n(z)-\RR_n(z)\, dz,
  \end{equation}
  where
  \begin{equation}
    \hat\RR_n(z) = \big(\Delta_n+(\RR_n(z))^{-1}\big)^{-1} = \RR_n(z)\big(\Delta_n\RR_n(z)+I\big)^{-1}.
  \end{equation}
  Since $\opnorm{\RR_n(z)}=1/\min_{1\le m\le d} |z-\lambda_m^{(n)}|<4/(3\eta)$ for $z\in Im(\Gamma)$ (see \cite[(5.3)]{HE15}), we have $\opnorm{\Delta_n\RR_n(z)}\le 4\opnorm{\Delta_n}/(3\eta)<2/3$. Hence (3.13) of \cite{HE15} yields
  $$ \hat\RR_n(z) = \RR_n(z) \Big(I+\sum_{m=1}^\infty (-\Delta_n\RR_n(z))^m\Big). $$
  Plug this into \eqref{eq:projdiff1} to conclude
  \begin{align}
    & \proj_{\hat W_j^{(n)}}- \proj_{W_j^{(n)}} \nonumber \\
     & = \frac 1{2\pi i} \oint_\Gamma \RR_n(z)\Delta_n \RR_n(z)\, dz -  \frac 1{2\pi i} \oint_\Gamma \RR_n(z) \sum_{m=2}^\infty (-\Delta_n\RR_n(z))^m \, dz.  \label{eq:projdiff2}
  \end{align}
  The operator norm of the second term on the right hand side is bounded by
  \begin{align}
   \frac{|\Gamma|}{2\pi} & \sup_{z\in Im(\Gamma)}\opnorm{\RR_n(z)} \sum_{m=2}^\infty \|\Delta_n\RR_n(z)\|_{op}^m  \nonumber\\
   & = \eta \frac{4}{3\eta} \Big(\frac{4\opnorm{\Delta_n}}{3\eta}\Big)^2 \frac{1}{1-2/3} \nonumber \\
   & \le \frac{64}{9\eta^2} \|\Delta_n\|_{op}^2 \nonumber \\
   & = O_P(k^{-1}), \label{eq:opnormbound}
  \end{align}
  where the last step follows from Theorem \ref{theo:empmixedmatconv}.
  Write $\proj_\ell^{(n)}$ as a shorthand for $\proj_{\Span{v_\ell^{(n)}}}$. By \cite[(5.2)]{HE15}, the leading term of the right hand side of \eqref{eq:projdiff2} equals
  \begin{align*}
    \sum_{\ell=1}^d & \sum_{m=1}^d \frac{1}{2\pi i} \oint_\Gamma \frac{1}{(\lambda_\ell^{(n)}-z)(\lambda_m^{(n)}-z)}\, dz \cdot \proj_\ell^{(n)}\Delta_n\proj_m^{(n)}.
  \end{align*}
  According to the residue theorem, the integral vanishes if both or neither of the eigenvalues $\lambda_\ell^{(n)}$ and $\lambda_m^{(n)}$ lie in the disk encircled by $\Gamma$, and it equals $(2\pi i)/(\lambda_m^{(n)}-\lambda_\ell^{(n)})=(2\pi i+o(1))/(\lambda_m^{(\infty)}-\mu_j)$ if $\lambda_\ell^{(n)}$ lies in the disk, but $\lambda_m^{(n)}$ does not. Combine this with \eqref{eq:projdiff2} and \eqref{eq:opnormbound} to obtain
  \begin{align*}
   &\proj_{\hat W_j^{(n)}}- \proj_{W_j^{(n)}}\\
    & = \sum_{\substack{\ell=1\\ \lambda_l=\mu_j}}^d \sum_{\substack{m=1\\ \lambda_m\ne\mu_j}}^d \frac{1}{\lambda_m^{(\infty)}-\mu_j}
  \big( \proj_\ell^{(n)}\Delta_n\proj_m^{(n)} + \proj_m^{(n)}\Delta_n\proj_\ell^{(n)}\big) +o_P(k^{-1/2})\\
  & = \sum_{\substack{i=1\\ i\ne j}}^{d_\infty} \frac{1}{\mu_i-\mu_j} \big( \proj_{W_j^{(n)}} \Delta_n
  \proj_{W_i^{(n)}} +  \proj_{W_i^{(n)}}\Delta_n  \proj_{W_j^{(n)}}\big)+o_P(k^{-1/2}).
  \end{align*}
   Now the assertion follows from Theorem \ref{theo:empmixedmatconv}, because  $\proj_{W_j^{(n)}} \to \proj_{W_j^{(\infty)}}$.
\end{proof}

\section{Asymptotic behavior of PCA projections and their excess risk}

\subsection{Unique optimal projection in the limit}

In this subsection we assume $\lambda_p^{(\infty)}>\lambda_{p+1}^{(\infty)}$, i.e.\ $p$ equals $N_{J_p}$ for some $J_p\in\{1,\ldots,d_\infty-1\}$ (see the discussion preceding Theorem \ref{theo:eigenproj} for the definition of $N_j$). Then
$$ V^*_\infty := V^*_{\infty,p} := \Span{\Big(\bigcup_{j=1}^{J_p} W_j^{(\infty)}\Big)}$$
is the unique $p$-dimensional subspace that maximizes
$$ M_\infty(\proj_V):=\E_\infty(\ang^T \proj_V\ang). $$

 Since its $i$th largest eigenvalue is a Lipschitz continuous function of a matrix, we conclude that $\lambda_p^{(n)}> \lambda_{p+1}^{(n)}$ for sufficiently large $n$ and  $V_{n,k}^*:=V_{n,k}^{*p}:=\Span\big(\bigcup_{j=1}^{J_p} W_j^{(n)}\big)$ is the unique maximizer of
$$ M_{n,k}(\proj_V) := \E_{n,k}(\ang^\top \proj_V \ang), \quad V\in\VV_p. $$
Likewise, with probability tending to 1, there is a unique maximizer  of
$$ \hat M_{n,k}(\proj_V) := \frac 1k \sum_{i=1}^n\ang_i^\top \proj_V \ang_i\Ind{R_i>\hat t_{n,k}}, \quad V\in\VV_p, $$
which by definition is the subspace obtained by our PCA procedure and equals
$$ \hat V_{n,k} := \hat V_{n,k}^p :=\Span\Big(\bigcup_{j=1}^{J_p} \hat W_j^{(n)}\Big). $$
Therefore, Theorem \ref{theo:eigenproj} implies a limit result for the difference between the PCA projection and the optimal projection.
\begin{cor} \label{cor:limitPCAproj}
  If $p=N_{J_p}$ for some $J_p\in\{1,\ldots,d_\infty-1\}$ and Condition (B) holds, then
  \begin{align*}
  k^{1/2} & \big(\proj_{\hat V_{n,k}} -\proj_{V_{n,k}^*} \big) \\
  & \to \sum_{j=1}^{J_p} \sum_{\ell=J_p+1}^{d_\infty} \frac 1{\mu_\ell-\mu_j} \big(\proj_{W_\ell^{(\infty)}} \Ub \proj_{W_j^{(\infty)}} + \proj_{W_j^{(\infty)}} \Ub \proj_{W_\ell^{(\infty)}}\big)=:\Wb_p
  \end{align*}
  weakly. The convergence holds jointly for all $p\in\{N_1,\ldots,N_{d_\infty-1}\}$.
\end{cor}
\begin{proof}
  Summation of the convergence result of Theorem \ref{theo:eigenproj} over $j\in\{1,\ldots,J_p\}$ yields weak convergence of $ k^{1/2} \big(\proj_{\hat V_{n,k}} -\proj_{V_{n,k}^*} \big)$ to
  \begin{align*}
     \sum_{j=1}^{J_p} & \sum_{\substack{\ell=1\\ \ell\ne j}}^{d_\infty} \frac 1{\mu_\ell-\mu_j} \big(\proj_{W_\ell^{(\infty)}} \Ub \proj_{W_j^{(\infty)}} + \proj_{W_j^{(\infty)}} \Ub \proj_{W_\ell^{(\infty)}}\big).
  \end{align*}
  Hence the assertion follows from
  \begin{align*}
   &\sum_{j=1}^{J_p}  \sum_{\substack{\ell=1\\ \ell\ne j}}^{J_p} \frac 1{\mu_\ell-\mu_j} \big(\proj_{W_\ell^{(\infty)}} \Ub \proj_{W_j^{(\infty)}} + \proj_{W_j^{(\infty)}} \Ub \proj_{W_\ell^{(\infty)}}\big)\\
  & = \sum_{1\le i<m\le J_p} \Big( \frac{1}{\mu_i-\mu_m}-\frac{1}{\mu_m-\mu_i}\Big)\big(\proj_{W_i^{(\infty)}} \Ub \proj_{W_m^{(\infty)}} + \proj_{W_m^{(\infty)}} \Ub \proj_{W_i^{(\infty)}}\big)\\
  & = 0.
  \end{align*}
\end{proof}

From Theorem \ref{theo:empmixedmatconv} we may also derive a limit distribution for the excess risk.
\begin{cor}  \label{cor:excessriskunique}
  If $p=N_{J_p}$  for some $J_p\in\{1,\ldots,d_\infty-1\}$ and Condition (B) holds, then
  \begin{align*}
    k & \big(M_{n,k}(\proj_{V_{n,k}^*})- M_{n,k}(\proj_{\hat V_{n,k}})\big) \to \sum_{j=1}^{J_p} \sum_{m=J_p+1}^{d_\infty} \frac{1}{\mu_j-\mu_m} \HSnorm{\proj_{W_j^{(\infty)}} \Ub \proj_{W_m^{(\infty)}}}
  \end{align*}
  weakly.
\end{cor}
\begin{proof}
  W.l.o.g.\ we may assume $\lambda_p^{(n)}>\hat\lambda_{p+1}^{(n)}$. Again we use $\proj_j^{(n)}$ as a shorthand for $\proj_{\Span v_j^{(n)}}$ and likewise write $\hat \proj_j^{(n)}$  for $\proj_{\Span \hat v_j^{(n)}}$. By Lemma 2.5 and Lemma 3.1 of \cite{ReissWahl20}, we have
  \begin{align}
    k & \big(M_{n,k}(\proj_{V_{n,k}^*})- M_{n,k}(\proj_{\hat V_{n,k}})\big) \nonumber\\
    & = \sum_{i=1}^p \sum_{\ell=p+1}^d \frac{\lambda_i^{(n)}}{(\lambda_i^{(n)}-\hat\lambda_\ell^{(n)})^2} \HSnorm{\proj_i^{(n)}\sqrt{k}\Delta_n \hat\proj_\ell^{(n)}} \nonumber\\
    & \hspace*{1cm} - \sum_{\ell=p+1}^d \sum_{i=1}^p \frac{\lambda_\ell^{(n)}}{(\lambda_\ell^{(n)}-\hat\lambda_i^{(n)})^2} \HSnorm{\proj_\ell^{(n)}\sqrt{k}\Delta_n \hat\proj_i^{(n)}}. \label{eq:excessconv}
  \end{align}
  From $\sqrt{k}\Delta_n\to \Ub$ weakly, $\lambda_i^{(n)}\to\mu_j$ for all $N_{j-1}<i\le N_j$, $\hat\lambda_\ell^{(n)}\to\mu_m$ in probability for all $N_{m-1}<\ell\le N_m$, as well as $\proj_{W_j^{(n)}}\to \proj_{W_j^{(\infty)}}$ and $\proj_{\hat W_j^{(n)}}\to \proj_{W_j^{(\infty)}}$ in probability for all $1\le j\le d_\infty$, we conclude that jointly for all $1\le j,m\le d_\infty$ with $j\ne m$
  \begin{align*}
  \sum_{i=N_{j-1}+1}^{N_j} & \sum_{\ell=N_{m-1}+1}^{N_m}  \frac{\lambda_i^{(n)}}{(\lambda_i^{(n)}-\hat\lambda_\ell^{(n)})^2} \HSnorm{\proj_i^{(n)}\sqrt{k}\Delta_n \hat\proj_\ell^{(n)}}\\
   & \to \frac{\mu_j}{(\mu_j-\mu_m)^2} \HSnorm{\proj_{W_j^{(\infty)}}\Ub\proj_{W_m^{(\infty)}}}
  \end{align*}
  weakly. Since $\big(\proj_{W_j^{(\infty)}}\Ub\proj_{W_m^{(\infty)}}\big)^\top = \proj_{W_m^{(\infty)}}\Ub\proj_{W_j^{(\infty)}}$, the right hand side of
  \eqref{eq:excessconv} converges weakly to
  \begin{align*}
      & \sum_{j=1}^{J_p} \sum_{m=J_p+1}^{d_\infty} \frac{\mu_j}{(\mu_j-\mu_m)^2} \HSnorm{\proj_{W_j^{(\infty)}}\Ub\proj_{W_m^{(\infty)}}} \\
      & \hspace{1cm} - \sum_{m=J_p+1}^{d_\infty} \sum_{j=1}^{J_p}  \frac{\mu_m}{(\mu_j-\mu_m)^2} \HSnorm{ \proj_{W_m^{(\infty)}}\Ub\proj_{W_j^{(\infty)}}} \\
      & = \sum_{j=1}^{J_p} \sum_{m=J_p+1}^{d_\infty} \frac{1}{\mu_j-\mu_m} \HSnorm{\proj_{W_j^{(\infty)}}\Ub\proj_{W_m^{(\infty)}}}.
  \end{align*}
\end{proof}
In particular, we see that the bound of the order $k^{-1/2}$ on the excess risk obtained in \cite{DS21} is very crude: if the optimal $p$-dimensional projection space is unique in the limit, then the bound decreases at the square root of the true rate.

\subsection{The general case}

We now drop the assumption that the maximizer of $M_\infty$ on $\VV_p$ is unique. Then, in general, one may not expect that the PCA projections converge, but one may hope for bounds on the excess risk. Reiss and Wahl \cite{ReissWahl20} analyzed the asymptotic behavior of the excess risk when all observations are drawn from a fixed distribution. In contrast, here the distribution of the exceedances varies with the sample size, which complicates the situation considerably, as a multiple eigenvalue in the limit may fan out quite arbitrarily for finite sample sizes.
\begin{cor} \label{cor:excessriskbound1}
   Define $J_p$ by the condition $N_{J_p-1}<p\le N_{J_p}$. Under Condition (B), one has for all maximizers $\hat V_{n,k}$ of $\hat M_{n,k}$ in $\VV_p$
  \begin{align}
       &\sup_{V\in\VV_p}M_{n,k}(\proj_{V_{n,k}^*})- M_{n,k}(\proj_{\hat V_{n,k}}) \nonumber\\
       & = O_P\bigg( \min_{1\le i\le p<j\le d+1} \Big( k^{-1}\Big(\frac{1}{\lambda_{i-1}^{(n)}-\lambda_{p+1}^{(n)}}+\frac{1}{\lambda_{p}^{(n)}-\lambda_{j}^{(n)}}\Big) +\lambda_i^{(n)}-\lambda_{j-1}^{(n)}\Big)\bigg) \label{eq:excessriskbound1}\\
       & = O_P\Big( k^{-1}+ \min\big(k^{-1/2},\lambda^{(n)}_{N_{J_p-1}+1}-\lambda^{(n)}_{N_{J_p}}\big)\Big), \label{eq:excessriskbound2}
        \end{align}
        with the convention $\lambda_0^{(n)}:=\infty$ and $\lambda_{d+1}^{(n)}:=-\infty$.
\end{cor}
 Note that the difference $\lambda^{(n)}_{N_{J_p-1}+1}-\lambda^{(n)}_{N_{J_p}}$ in \eqref{eq:excessriskbound2} can be bounded by \linebreak $2\hsnorm{\Sigma_{n,k}-\Sigma_\infty}$. Hence we again obtain a rate of convergence $k^{-1}$ if we assume $\hsnorm{\Sigma_{n,k}-\Sigma_\infty}=O(k^{-1})$. More generally, the excess risk converges fast if there is a cluster of eigenvalues around $\lambda_p^{(n)}$ that is closely concentrated and clearly separated from all other eigenvalues. In the worst case, the eigenvalues corresponding to $\lambda_p$ in the limit widely fan out for finite sample sizes, which may lead to a slow rate of $k^{-1/2}$ of the upper bounds.
\begin{proof}
  By Theorem \ref{theo:empmixedmatconv}, to each $\eps>0$ there exists $C>0$ such that
  $$ \P\Big\{ \max_{1\le i\le d} |\hat \lambda^{(n)}_i- \lambda^{(n)}_i|>C k^{-1/2} \text{ or } \hsnorm{\Delta_n}>Ck^{-1/2}\Big\}\le\eps
  $$
  for $n$ sufficiently large.

  Let $(i_n,j_n)$ be a pair for which the minimum in \eqref{eq:excessriskbound1} is attained. W.l.o.g.\ we may assume $\lambda_{i_n-1}^{(n)}-\lambda_{p+1}^{(n)}>2Ck^{-1/2}$. If this is not the case, then replacing $i_n$ with $\tilde i_n := \max\big\{i<i_n\mid \lambda_{i-1}^{(n)}-\lambda_{p+1}^{(n)}>2Ck^{-1/2}\big\}$ will decrease the term $k^{-1}/(\lambda_{i-1}^{(n)}-\lambda_{p+1}^{(n)})$ and increase the term $\lambda_{i}^{(n)}-\lambda_{j-1}^{(n)}$ by at most $2Ck^{-1/2}=O\big(k^{-1}/(\lambda_{i-1}^{(n)}-\lambda_{p+1}^{(n)})\big)$; it will thus not increase the order of \eqref{eq:excessriskbound1}. Likewise, we may assume $\lambda_p^{(n)}-\lambda_{j_n}^{(n)}>2Ck^{-1/2}$.

  Therefore, with probability of at least $1-\eps$, we have for all $i<i_n$ and $j\ge p+1$
  $$ (\lambda_i^{(n)}-\hat\lambda_j^{(n)})^2 \ge  (\lambda_i^{(n)}-\lambda_{p+1}^{(n)}-Ck^{-1/2})^2\ge \frac 14  (\lambda_i^{(n)}- \lambda_{p+1}^{(n)})^2,
  $$
  and for all $j\ge j_n$ and $i\le p$
  $$ (\hat\lambda_i^{(n)}-\lambda_j^{(n)})^2 \ge  (\lambda_p^{(n)}-\lambda_{j}^{(n)}-Ck^{-1/2})^2\ge \frac 14  (\lambda_p^{(n)}- \lambda_{j}^{(n)})^2.
  $$
  Hence, with the notation introduced in the proof of Corollary \ref{cor:excessriskunique}, by combining Lemma 2.5 and Lemma 3.1 of \cite{ReissWahl20} we obtain with probability of at least $1-\eps$
  \begin{align*}
       &\sup_{V\in\VV_p}M_{n,k}(\proj_{V_{n,k}^*})- M_{n,k}(\proj_{\hat V_{n,k}})\\
       & \le \sum_{i=1}^{i_n-1} (\lambda_i^{(n)}-\lambda_p^{(n)}) \sum_{j=p+1}^d \frac{\HSnorm{\proj_i^{(n)}\Delta_n\hat\proj_j^{(n)}}}{(\lambda_i^{(n)}-\hat\lambda_j^{(n)})^2}\\
       & \quad +
       \sum_{j=j_n}^{d} (\lambda_p^{(n)}-\lambda_j^{(n)}) \sum_{i=1}^p \frac{\HSnorm{\proj_j^{(n)}\Delta_n\hat\proj_i^{(n)}}}{(\hat\lambda_i^{(n)}-\lambda_j^{(n)})^2}
        + \sum_{i=i_n}^{j_n-1} |\lambda_i^{(n)}-\lambda_p^{(n)}|\\
       & \le 4C^2 k^{-1} \bigg[ (d-p) \sum_{i=1}^{i_n-1} \frac{\lambda_i^{(n)}-\lambda_p^{(n)}}{(\lambda_i^{(n)}- \lambda_{p+1}^{(n)})^2} + p  \sum_{j=j_n}^{d} \frac{\lambda_p^{(n)}-\lambda_j^{(n)}}{(\lambda_p^{(n)}- \lambda_{j}^{(n)})^2}\bigg] \\
       & \quad + (j_n-i_n)(\lambda_{i_n}^{(n)}-\lambda_{j_n-1}^{(n)})\\
       & \le 4C^2 k^{-1} \bigg[ (d-p) \frac{i_n-1}{\lambda_{i_n-1}^{(n)}-\lambda_{p+1}^{(n)}}+p \frac{d-j_n+1}{\lambda_p^{(n)}-\lambda_{j_n}^{(n)}}\bigg]+ (j_n-i_n)(\lambda_{i_n}^{(n)}-\lambda_{j_n-1}^{(n)}),
  \end{align*}
  which proves the bound \eqref{eq:excessriskbound1}.

  If $\lambda^{(n)}_{N_{J_p-1}+1}-\lambda^{(n)}_{N_{J_p}}\le k^{-1/2}$, then the second bound is obtained by choosing $i=N_{J_p-1}+1$ and $j=N_{J_p}+1$. Else choose $i=\max\{l\le p\mid \lambda_l^{(n)}-\lambda_{p+1}^{(n)}> k^{-1/2}\}+1$ and $j=\min\{l>p \mid \lambda_{p}^{(n)}-\lambda_l^{(n)}> k^{-1/2}\}$. (Note that by assumption both sets are not empty.) Then both denominators in \eqref{eq:excessriskbound1}
  are greater than $k^{-1/2}$ and $\lambda_i^{(n)}-\lambda_{j-1}^{(n)}$ is at most $2k^{-1/2}$.
\end{proof}


\section{Local empirical processes}
\label{sect:locempproc}

Throughout this section, we assume that $p=N_{J_p}$ for some $J_p\in\{1,\ldots, d_\infty-1\}$, so that in the limit the optimal projection is unique. According to Corollary \ref{cor:excessriskunique}, the PCA projection differs from the optimal projection $\proj_{V_{n,k}^*}$ only by a matrix of order $k^{-1/2}$. It is thus natural to examine the asymptotic behavior of $\hat M_{n,k}$ on a neighborhood of $\proj_{V_{n,k}^*}$ that shrinks with the rate $k^{-1/2}$.

Note that every projection on a subspace $V\in\VV_p$ is of the form $S^\top\proj_{V_{n,k}^*}S$ for some rotation matrix $S\in SO(n)$. It is well known that the special orthogonal group consists of all matrix exponentials of skewed symmetric matrices. Let
$$ \MM_S:=\{A\in\R^{d\times d}\mid A^\top=-A\}. $$
Then for all $V\in\VV_p$ there exists $A\in \MM_S$ such that
$ \proj_V = \e^{-A}\proj_{V_{n,k}^*}\e^A $
and any $A\in\MM_S$ defines a projection onto some $V\in\VV_p$. Hence, it suggests itself to use the following local parametrization:
$$ \tilde\proj_{n,A} := \e^{-k^{-1/2}A}\proj_{V_{n,k}^*}\e^{k^{-1/2}A}. $$
In particular $\tilde\proj_{n,0}=\proj_{V_{n,k}^*}$. A Taylor expansion of the matrix exponential yields
\begin{align}
  \tilde\proj_{n,A}
  & = \Big(I-k^{-1/2}A+\frac 12 k^{-1} A^2\Big) \proj_{V_{n,k}^*} \Big(I+k^{-1/2}A+\frac 12 k^{-1} A^2\Big) +O(k^{-3/2}) \nonumber\\
  & = \proj_{V_{n,k}^*} + k^{-1/2} (\proj_{V_{n,k}^*} A-A\proj_{V_{n,k}^*})  \nonumber \\
   & \hspace{1.5cm} +k^{-1}\Big(\frac{\proj_{V_{n,k}^*}A^2+A^2\proj_{V_{n,k}^*}}{2}- A\proj_{V_{n,k}^*}A\Big)+O(k^{-3/2}),
     \label{eq:projexpan}
\end{align}
(with $I$ denoting the identity matrix) which holds uniformly over compact sets in $\MM_S$. Note, however, that different matrices $A$ may lead to the same or asymptotically equivalent projections (up to terms of smaller order than $k^{-1/2}$, which will be asymptotically indistinguishable). More precisely, for any two matrices $A,B\in\MM_S$ one has
\begin{align*}
  & \hsnorm{\tilde\proj_{n,A}-\tilde\proj_{n,B}}=o(k^{-1/2})\\
   \iff & \hsnorm{\proj_{V_{n,k}^*} (A-B)-(A-B)\proj_{V_{n,k}^*}}=o(1) \\
   \iff & \proj_{V_\infty^*} (A-B)=(A-B)\proj_{V_\infty^*} \\
   \iff & (A-B)V_\infty^*\subset V_\infty^* \text{ and } (A-B)(V_\infty^*)^\perp\subset (V_\infty^*)^\perp.
\end{align*}
In what follows, we will thus parametrize the projections with the set of matrices
\begin{align*}
\tilde\MM_S &:= \big\{A\in\R^{d\times d} \mid A^\top=-A, AV_\infty^*\subset (V_\infty^*)^\perp,A(V_\infty^*)^\perp\subset V_\infty^*\big\}.
\end{align*}
Then, for all $B\in\MM_S$, there exists $A\in\tilde\MM_S$  such that $\hsnorm{\tilde\proj_{n,A}-\tilde\proj_{n,B}}=o(k^{-1/2})$ (namely $A=\proj_{(V_\infty^*)^\perp} B \proj_{V_\infty^*} + \proj_{V_\infty^*} B \proj_{(V_\infty^*)^\perp}$). It is the main aim of this section to describe the asymptotic behavior of $M_{n,k}$ and $\hat M_{n,k}$  locally in the neighborhood of $\proj_{V_{n,k}^*}$ and $\proj_{\hat V_{n,k}}$ using such a local parametrization.

For notational simplicity, in what follows we drop the index $\infty$ when dealing with the limit model. In particular, $\lambda_i$ denotes the $i$th largest eigenvalue of $\Sigma=\Sigma_\infty$ and $V^*=V_\infty^*\in\VV_p$ denotes the optimal projection space in the limit model. Moreover, we use $V^\perp$, $\proj^*$, $\proj^\perp$ and $\proj_i$ ($1\le i\le p$) as short hands for $(V^*_\infty)^\perp$, $\proj_{V_\infty^*}$, $\proj_{(V_\infty^*)^\perp}$ and $\proj_{\Span(v_i^{(\infty)})}$, respectively. Define linear mappings
\begin{align*}
  S_\lambda:\; & \tilde\MM_S\to \MM^* :=\{A\in\R^{d\times d}\mid AV^*=\{0\}, AV^\perp\subset V^*\}\\
    & S_\lambda(A) := \sum_{i=1}^p \sum_{j=p+1}^d \sqrt{\lambda_i-\lambda_j} \proj_iA\proj_j,\\[0.5ex]
  T_\lambda:\; & \R^{d\times d}\to\MM^*\\
    & T_\lambda(B) := \sum_{i=1}^p \sum_{j=p+1}^d \frac{\proj_iB\proj_j}{\sqrt{\lambda_i-\lambda_j}},\\[0.5ex]
  \bar T_\lambda: \;& \MM^*\to\tilde\MM_S; \quad \bar T_\lambda(B) := T_\lambda(B)-(T_\lambda(B))^\top.
\end{align*}
Direct calculations show that $\bar T_\lambda$ is the inverse mapping to $S_\lambda$.

The following result describes the asymptotic behavior of the risk and the empirical risk of projections in the neighborhood of the optimal projection.
\begin{theorem} \label{theo:localproc1}
  If Condition (B) holds and there is a unique optimal projection in the limit, then
  \begin{align}
    k&\big(M_{n,k}(\tilde\proj_{n,A})-M_{n,k}(\proj_{V_{n,k}^*})\big)\;\longrightarrow
    \; \hsprod{\Sigma}{A^2(\proj^*-\proj^\perp)}  \label{eq:localrisk}\\
   k&\big(\hat M_{n,k}(\tilde\proj_{n,A})-\hat M_{n,k}(\proj_{V_{n,k}^*})\big) \nonumber\\
   &\;\longrightarrow
    \;2\hsprod{\Ub}{\proj^* A} + \hsprod{\Sigma}{A^2(\proj^*-\proj^\perp)} \quad  \text{weakly,} \label{eq:localemprisk}
     \end{align}
     uniformly over compact sets of matrices $A\in\tilde\MM_S$.
  The terms arising in the limit can be expressed in terms of $S_\lambda$ and $T_\lambda$ as follows:
  \begin{align*}
     \hsprod{\Sigma}{A^2(\proj^*-\proj^\perp)} & = - \HSnorm{S_\lambda(A)} \\
     \hsprod{\Ub}{\proj^* A} & = \hsprod{T_\lambda(\Ub)}{S_\lambda(A)}.
  \end{align*}
 The unique maximizer of the limit process in \eqref{eq:localemprisk} is $A^*=\bar T_\lambda(T_\lambda(\Ub))$ with the corresponding maximum being equal to $\HSnorm{T_\lambda(\Ub)}$, and the corresponding limit in \eqref{eq:localrisk} equals $-\HSnorm{T_\lambda(\Ub)}$.
\end{theorem}
\begin{proof}
  To establish \eqref{eq:localrisk}, we express the difference as a product in the space of Hilbert-Schmidt operators, using \eqref{eq:projexpan}:
  \begin{align*}
    k &\big(M_{n,k}(\tilde\proj_{n,A})-M_{n,k}(\proj_{V_{n,k}^*})\big)\\
    & = k\hsprod{\Sigma_{n,k}}{\tilde\proj_{n,A}-\proj_{V_{n,k}^*}} \\
    & = k^{1/2} \hsprod{\Sigma_{n,k}}{\proj_{V_{n,k}^*}A-A\proj_{V_{n,k}^*}} \\
     & \quad + \hsprod{\Sigma_{n,k}}{\frac{\proj_{V_{n,k}^*}A^2+A^2\proj_{V_{n,k}^*}}{2}-A\proj_{V_{n,k}^*}A} +O(k^{-1/2}).
  \end{align*}
   Because $\Sigma_{n,k}$ and $ \proj_{V_{n,k}^*}=\proj_{V_{n,k}^*}^\top$ commute, the first term on the right hand side  vanishes:
   \begin{align*}
     \hsprod{\Sigma_{n,k}}{\proj_{V_{n,k}^*}A} & = \hsprod{\proj_{V_{n,k}^*}^\top\Sigma_{n,k}}{A} = \hsprod{\Sigma_{n,k}\proj_{V_{n,k}^*}^\top}{A} \\
     & = \hsprod{\Sigma_{n,k}}{A\proj_{V_{n,k}^*}}.
   \end{align*}
   Since $\Sigma_{n,k}\to\Sigma$ and $\proj_{V_{n,k}^*}\to\proj^*$, convergence \eqref{eq:localrisk} follows if we show that
   $$ \hsprod{\Sigma}{\frac{\proj^*A^2+A^2\proj^*}{2}-A\proj^*A} = \hsprod{\Sigma}{A^2(\proj^*-\proj^\perp)}. $$
   To this end, note that for all symmetric matrices $S\in\R^{d\times d}$
   $$ \hsprod{S}{\proj^* A^2} = \hsprod{S^\top}{(\proj^*A^2)^\top} = \hsprod{S}{A^2\proj^*} $$
   and $\proj^*A=A\proj^\perp$, by the definition of $\tilde \MM_S$, which gives the desired equality.

   To prove \eqref{eq:localemprisk}, recall the notation $\Delta_n:=\hat\Sigma_{n,k}-\Sigma_{n,k}$. In view of \eqref{eq:projexpan} and Theorem \ref{theo:empmixedmatconv}, we have
   \begin{align*}
     k &  \big(\hat M_{n,k}(\tilde\proj_{n,A})-\hat M_{n,k}(\proj_{V_{n,k}^*}) \big)- k \big(M_{n,k}(\tilde\proj_{n,A})-M_{n,k}(\proj_{V_{n,k}^*})\big) \nonumber\\
       & = \hsprod{k^{1/2}(\hat\Sigma_{n,k}-\Sigma_{n,k})}{k^{1/2}(\tilde\proj_{n,A}-\proj_{V_{n,k}^*})} \\
       & = \hsprod{k^{1/2}\Delta_n}{\proj_{V_{n,k}^*}A-A\proj_{V_{n,k}^*}}+O_P(k^{-1/2}). 
   \end{align*}

   Since $\proj_{V_{n,k}^*}\to\proj^*$ and $\hsprod{S}{\proj^* A} = \hsprod{S^\top}{(\proj^*A)^\top} = \hsprod{S}{-A\proj^*} $ for all  symmetric matrices $S\in\R^{d\times d}$, Theorem \ref{theo:empmixedmatconv}, the continuous mapping theorem  and  \eqref{eq:localrisk} imply \eqref{eq:localemprisk}.

   It remains to express the limits in terms of $S_\lambda$ and $T_\lambda$.
   To this end, check that, for all $A\in\tilde\MM_S$, $i,\ell\in\{1,\ldots,p\}$ and $j,m\in\{p+1,\ldots,d\}$
    \begin{align*}
    \hsprod{\proj_i\Ub\proj_j}{\proj_\ell A\proj_m}
    & = \tr\big(\proj_i\Ub\proj_j(\proj_\ell A\proj_m)^\top\big)\\
     & = \tr(\Ub\proj_jA^\top\proj_\ell\proj_i)\Ind{j=m} \\
     & = \hsprod{\Ub}{\proj_iA\proj_j}\Ind{i=\ell,j=m},
   \end{align*}
   and thus
   \begin{align}
     \hsprod{T_\lambda(\Ub)}{S_\lambda(A)}
     & = \sum_{i=1}^p\sum_{j=p+1}^d \sum_{\ell=1}^p \sum_{m=p+1}^d \frac{\sqrt{\lambda_\ell-\lambda_m}}{\sqrt{\lambda_i-\lambda_j}} \hsprod{\proj_iU\proj_j}{\proj_\ell A\proj_m}  \nonumber \\
     & = \sum_{i=1}^p\sum_{j=p+1}^d \hsprod{\Ub}{\proj_iA\proj_j}  \nonumber \\
     & = \hsprod{\Ub}{\proj^*A\proj^\perp} \nonumber \\
     & = \hsprod{\Ub}{\proj^*A}. \label{eq:loclimitrandom}
   \end{align}
   In the last step again the definition of $\tilde\MM_S$ has been used. Similarly, we conclude from
   \begin{align*}
     \hsprod{\proj_i A\proj_j}{\proj_\ell A\proj_m} & = \tr(\proj_jA^\top\proj_i A) \Ind{i=\ell,j=m}=-\tr(\proj_jA\proj_i A) \Ind{i=\ell,j=m}\\
     & =-\tr(\proj_iA\proj_j A) \Ind{i=\ell,j=m}
   \end{align*}
    and $\proj_j A \proj_i=0$ if $i,j\in\{1,\ldots,p\}$ or $i,j\in\{p+1,\ldots, d\}$ that
   \begin{align*}
    \HSnorm{S_\lambda(A)} & = \sum_{i=1}^p\sum_{j=p+1}^d \sum_{\ell=1}^p \sum_{m=p+1}^d  \sqrt{\lambda_i-\lambda_j} \sqrt{\lambda_\ell-\lambda_m}\hsprod{\proj_iA\proj_j}{\proj_\ell A\proj_m} \\
    & = -\sum_{i=1}^p\sum_{j=p+1}^d (\lambda_i-\lambda_j)\tr(\proj_jA\proj_iA)\\
    & = -\sum_{i=1}^d \lambda_i \sum_{j=1}^d \tr(\proj_jA\proj_iA)\big(\Ind{i\le p}-\Ind{i> p}\big)\\
    & = -\sum_{i=1}^d \lambda_i \tr(A^2\proj_i)\big(\Ind{i\le p}-\Ind{i> p}\big)\\
    & = -\sum_{i=1}^d \lambda_i \tr\big(A^2(\proj^*-\proj^\perp)\proj_i\big)\\
    & = -\sum_{i=1}^d \lambda_i \hsprod{\proj_i}{A^2(\proj^*-\proj^\perp)}\\
    & = -\hsprod{\Sigma}{A^2(\proj^*-\proj^\perp)}.
    \end{align*}

    Obviously, any maximizer $A$ of the limit process $2\hsprod{T_\lambda(\Ub)}{S_\lambda(A)}$ \linebreak $-\HSnorm{S_\lambda(A)}$ must fulfill $S_\lambda(A)=cT_\lambda(\Ub)$ for some (possibly random) $c\in\R$. But then the limit equals
    $2c\HSnorm{T_\lambda(\Ub)}-c^2\HSnorm{T_\lambda(\Ub)}$, which is maximal for $c=1$. Because $\bar T_\lambda$ is the inverse to $S_\lambda$, we finally obtain $A^*=\bar T_\lambda(T_\lambda(\Ub))$ as the unique maximizer of the limit process, leading to the maximum value $\HSnorm{T_\lambda(\Ub)}$ in \eqref{eq:localemprisk} and to $-\HSnorm{T_\lambda(\Ub)}$ in \eqref{eq:localrisk}.
\end{proof}
\begin{remark}
  It is easily seen that the maximum of the limit process equals the limit in Corollary \ref{cor:excessriskunique} and that the unique maximizer can be expressed as
  $$ A^* = \sum_{i=1}^p \sum_{j=p+1}^d \frac{\proj_i U \proj_j-\proj_j U\proj_i}{\lambda_i-\lambda_j}.
  $$
  Therefore,
  \begin{align*}
    k^{1/2} \big( \tilde\proj_{n,A^*}-\proj_{V_{n,k}^*}\big)
      & \; \longrightarrow\; \proj^*A^*-A^*\proj^* \; = \sum_{i=1}^p \sum_{j=p+1}^d \frac{\proj_i U \proj_j+\proj_j U\proj_i}{\lambda_i-\lambda_j},
  \end{align*}
  which equals the limit in Corollary \ref{cor:limitPCAproj}. These results could hence also be concluded from Theorem  \ref{theo:localproc1}, by invoking a suitable argmax-theorem, provided one first proves tightness of $k^{1/2}\big(\proj_{\hat V_{n,k}}-\proj_{V_{n,k}^*}\big)$.
\end{remark}

While Theorem \ref{theo:localproc1} gives a neat description of the asymptotic behavior of projections in a neighborhood of the optimal projection, it is difficult to appy to any given projection matrix, because one would need to know $\proj_{V_{n,k}^*}$ in order to calculate the pertaining local parameter $A$. To avoid this problem, one might also think of a different local parametrization centered at the known (random) PCA projection:
$$ \hat\proj_{n,B} := \e^{-k^{-1/2}B}\proj_{\hat V_{n,k}}\e^{k^{-1/2}B}, \quad B\in\tilde \MM_S. $$
Arguments completely analogous to the ones employed in the proof of Theorem \ref{theo:localproc1}, but with the roles of $(\Sigma_{n,k},V_{n,k})$ and $(\hat\Sigma_{n,k},\hat V_{n,k})$ interchanged, yield the following result.
\begin{theorem} \label{theo:localproc2}
  If Condition (B) is met and there is a unique optimal projection in the limit, then
  \begin{align}
    k&\big(\hat M_{n,k}(\hat\proj_{n,B})-\hat M_{n,k}(\proj_{\hat V_{n,k}})\big)\;\longrightarrow
    \; \hsprod{\Sigma}{B^2(\proj^*-\proj^\perp)}  \label{eq:localrisk2}\\
   k&\big(M_{n,k}(\hat\proj_{n,B})-M_{n,k}(\proj_{\hat V_{n,k}})\big) \nonumber\\
   &\;\longrightarrow
    \; -2\hsprod{\Ub}{\proj^* B} + \hsprod{\Sigma}{B^2(\proj^*-\proj^\perp)} \quad  \text{weakly,} \label{eq:localemprisk2}
     \end{align}
     uniformly over compact sets of matrices $B\in\tilde\MM_S$.
\end{theorem}

\section{Choice of the dimension}

So far we have assumed that the dimension $p$ of the projection subspace is given. In this section we discuss an approach to select $p$, which is then assessed in a simulation study. For ease of presentation, in this chapter we assume that the observations are re-scaled to the unit sphere, that is $\ang_i=X_i/\|X_i\|$, $1\le i\le n$.

\subsection{Selection procedure}

If the angular measure is concentrated on a proper subspace of $\R^d$ (intersected with the unit sphere), then a natural goal is to estimate the dimension of the smallest subspace supporting $\ang$ in the limit model. Note that, unlike in the classical setting, where the whole underlying distribution is to be approximated and hence all observations lie in any supporting subspace, figuring out the smallest $p$ such that the angular measure is concentrated on some $V\in\VV_p$ is a non-trivial task, because the observations are not drawn from the limit model we want to analyze.

More generally, if the support of the angular measure may span the whole $\R^d$, then one might aim at finding the smallest $p$ such that the true reconstruction error of the projection on some optimal subspace $V\in\VV_p$ lies below a given threshold. Equivalently, one may fix some $\tau\in (0,1)$ and try to estimate $p_\tau := \min\big\{p\in\{1,\ldots,d\}\mid M_\infty(\proj_V):=\E_\infty(\ang^\top\proj_V\ang)\ge\tau$ for some $V\in\VV_p\big\}$. Since $\sup_{V\in\VV_p}M_\infty(\proj_V)$ is unknown and must be replaced with $\hat M_{n,k}(\proj_{\hat V_{n,k}})$, the estimation error has to be taken into account. Without further model assumptions, it is impossible to say anything about the bias part which arises from the fact that the function $M_\infty$ in the limit model differs from $M_{n,k}$ which can be statistically analyzed. In what follows, we always assume that this bias is asymptotically negligible, which implies that $k$ has been chosen sufficiently small.

Assuming the existence of a unique maximizer  in the limit and thus eventually also a unique maximizer $V=V_{n,k}^{*p}$ of dimension $p$ of $M_{n,k}(\proj_V)$ and (with probability arbitrarily close to 1) a unique maximizer $V=\hat V_{n,k}^p$ of $\hat M_{n,k}(\proj_V)$, the remaining estimation error equals
\begin{align*}
  \hat M_{n,k}& (\proj_{\hat V_{n,k}^p})-M_{n,k}(\proj_{V_{n,k}^{*p}}) \\
    & = \big(\hat M_{n,k} (\proj_{\hat V_{n,k}^p})- M_{n,k} (\proj_{\hat V_{n,k}^p})\big) + \big(M_{n,k} (\proj_{\hat V_{n,k}^p})-M_{n,k}(\proj_{V_{n,k}^{*p}})\big).
\end{align*}
The second term on the right-hand side is the excess risk, which is of the order $k^{-1}$. In contrast, the first term, which can be rewritten as $\hsprod{\Delta_n}{\proj_{\hat V_{n,k}^p}}$ is usually of the order $k^{-1/2}$. More precisely, in view of our results from Sections \ref{sect:convempmom} and  \ref{sect:locempproc}, we see that
\begin{align*}
  k^{1/2} & \big(\hat M_{n,k} (\proj_{\hat V_{n,k}^p})- M_{n,k} (\proj_{\hat V_{n,k}^p})\big)
    \to \hsprod{U}{\proj^{(*p)}} = \sum_{i=1}^p v_i^\top U v_i,
\end{align*}
where $\proj^{(*p)}$ denotes the optimal projection in the limit onto a $p$ dimensional subspace.
The limit random variable is centered normal  with variance
\begin{align*}
 \sigma_p^2 & = \sum_{i=1}^p \sum_{j=1}^p Cov_\infty(v_i^\top \ang\ang^\top v_i,v_j^\top \ang\ang^\top v_j) \\
   & = Var_\infty\Big(\sum_{i=1}^p (v_i^\top\ang)^2\Big)\\
   & = Var_\infty(\|\proj^{(*p)}\ang\|^2).
\end{align*}
A natural estimator for $\sigma_p^2$ is the empirical variance of the squared norm of  $\proj_{\hat V_{n,k}^{p}}\ang$, that is
$$ \hat\sigma_p^2 := \frac{1}{k-1} \sum_{l=1}^n \Big( \| \proj_{\hat V_{n,k}^{p}}\ang_l\|^2\Ind{R_l>\hat t_{n,k}}
-\sum_{i=1}^p\hat \lambda_i^{(n)}\Big)^2.
$$
So if $\hat M_{n,k}(\proj_{\hat V_{n,k}^{p}})=\sum_{i=1}^p \hat \lambda_i^{(n)}$ exceeds $\tau+k^{-1/2}z_\beta\hat\sigma_p$ (with $z_\beta$ denoting the standard normal $\beta$-quantile), then one can be sure with approximate probability $\beta$, that the true reconstruction error is at most $1-\tau$.

To sum up, we propose to fix some $\tau,\beta\in (0,1)$ (typically close to 1) and choose the dimension as
\begin{equation} \label{eq:choice_p}
 \hat p := \min\Big\{p\le d \mid \sum_{i=1}^p \hat \lambda_i^{(n)}>\tau+k^{-1/2}z_\beta\hat\sigma_p\Big\}.
\end{equation}
Strictly speaking, this definition implicitly invokes a multiple testing problem, which partly invalidates the above interpretation of the tuning parameter $\beta$, but other sources of distortion like the aforementioned bias seem more relevant.

\subsection{Simulation study}

In what follows we examine the performance of this approach in a small simulation study. We largely follow the setting introduced in \cite{DS21}. More concretely, we sample from three different dependence models:
\begin{itemize}
  \item the first $p$ coordinates of the random vectors have a Gumbel copula (and the remaining coordinates equal 0);
  \item the first $p$ coordinates are drawn from a Dirichlet model;
  \item the first $p$ coordinates are drawn from a Dirichlet model and then the vector is randomly rotated by an angle that is uniformly distributed on $[-\pi/10,\pi/10]$ in the plane spanned by two unit vectors randomly chosen among the first $p$ and the last $d-p$ unit vectors, respectively.
\end{itemize}
In addition, the modulus of some normally distributed noise is added, so that the observations are not drawn from the limit model. In all models, the marginals are Fr\'{e}chet distributed with some tail index $\alpha\in\{1,2\}$. We examine models of moderate dimension $d=10$ and $p=2$ and high dimensional models with $d=100$ and $p=5$. Note that in the first two models the angular measure is concentrated on a $p$ dimensional subspace, whereas in the last model the mass is spread out over a neighborhood of such a subspace. See \cite{DS21} for a detailed description of these models.

In each setting, we draw $1\,000$ samples  of size $n=1\,000$. For each sample, we perform the PCA procedure with the ``true'' dimension $p$ as well as with dimension $\hat p$ selected according to \eqref{eq:choice_p}, with $\tau=\beta=0.95$, using different values $k$ of the number of observations considered extreme, which results in estimators  of the angular measure
\begin{equation} \label{eq:hatHnkPCA}
 \hat H_{n,k}^{PCA} := \frac 1k \sum_{i=1}^n \delta_{\angf_{t_{n,k}}(\proj_{\hat V_{n,k}^{\tilde p}} X_i)},
\end{equation}
with $\tilde p\in\{p,\hat p\}$. In \cite{DS21} it is shown that estimators of the angular measure using PCA often perform better if the selection of the lower dimensional subspace is based on a smaller number $\tilde k$ of extreme observations than the number used by the empirical angular measure. This can be explained by the fact that the most extreme observations are least perturbed by the light tailed noise and lie thus closest to the lower dimensional subspace. However, a very small number of extreme observations does not contain sufficient information to estimate the angular measure accurately. Thus, in addition, we consider the PCA projection $\proj_{\hat V_{n,\tilde k}^{\tilde p}}$ and the dimension selector $\hat p$ based on a small number $\tilde k$ of largest observations, but then calculate the estimate of the angular measure using a (potentially) larger number $k$ of observations to obtain
\begin{equation} \label{eq:hatHnkktildePCA}
 \hat H_{n,k,\tilde k}^{PCA}:= \frac 1k \sum_{i=1}^n \delta_{\angf_{t_{n,k}}(\proj_{\hat V_{n,\tilde k}^{\tilde p}} X_i)}.
\end{equation}

In most applications, estimating the angular measure is not the main goal of the statistical analysis, but it is just an important step in order to estimate some other quantities. To evaluate the performance of the proposed PCA procedure, for $X=(X^1,\ldots,X^d)$, we thus calculate resulting estimators of the following four probabilities in the limit model:
\begin{itemize}
   \item[(i)] $\lim_{u\to\infty} \P(p^{-1}\sum_{1\le j\le p} X^j/\|X\|>t_{(i)} \mid \|X\|>u)=H\{x\mid p^{-1} \sum_{j=1}^p x^j>t_{(i)}\}$ for some $t_{(i)}\in (0,p^{-1/2})$
   \item[(ii)] $\lim_{u\to\infty} \P(\min_{1\le j\le p}X^j>u,\max_{p+1\le j\le d} X^j\le u \mid \|X\|>u)=$\\ \hspace*{2cm} $ \int \big((\min_{1\le j\le p}x^j)^\alpha-(\max_{p+1\le j\le d}x^j)^\alpha\big)^+\, H(dx)$
    \item[(iii)] $ \lim_{u\to\infty} \P(X^1>u\mid \max_{1\le j\le d}X^j>u) = $ \\
     \hspace*{2cm} $\int (x^1)^\alpha\, H(dx)/\int (\max_{1\le j\le p}x^j)^\alpha\, H(dx)$
   \item[(iv)] $ \lim_{u\to\infty} \P(\min_{1\le j\le d}X^j>u\mid \|X\|>u) = \int (\min_{1\le j\le d}x^j)^\alpha\, H(dx)$
\end{itemize}
In \cite{DS21} a detailed discussion of the interpretation of these parameters is given.

We will first discuss the settings with moderate dimension ($d=10$, $p=2$) and then the high dimensional models ($d=100$, $p=5$).

\subsubsection*{Moderate dimensional models}

\begin{figure}
\includegraphics[width=0.9\textwidth]{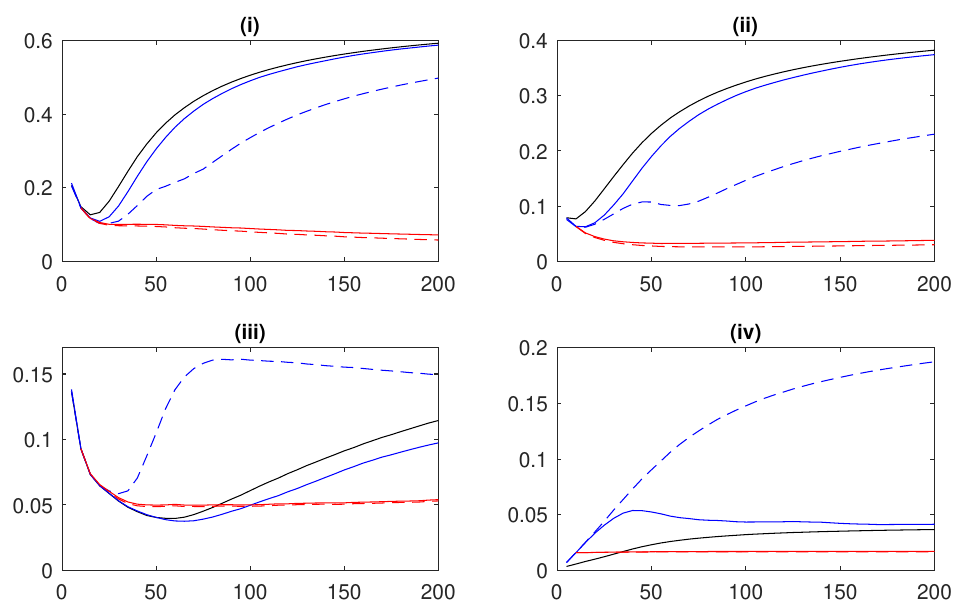}
\caption{RMSE of the estimators of the probabilities (i)--(iv) based on $\hat H_{n,k}$ (black, solid), $\hat H_{n,k}^{PCA}$ (blue, dashed) and $\hat H_{n,k,10}^{PCA}$ (red, dashed) versus $k$ in the Dirichlet model with fixed dimensions $p=2$ and $d=10$. The colored solid lines indicate the corresponding RMSE when the dimension is chosen as in \eqref{eq:choice_p}.}
 \label{fig:RMSEdiri3_10_2}
\end{figure}

We start with the Dirichlet model with all Dirichlet parameters $\alpha_i$, $1\le i\le p$, equal to 3 and unit Fr\'{e}chet margins. Figure \ref{fig:RMSEdiri3_10_2} shows the root mean squared error (RMSE) of the estimators of the parameters (i)--(iv) (with $t_{(i)}=0.65$)  as a function of $k$. The true values of the parameters are (approximately) $0.6838$, $0.4558$, $0.7619$ and $0$, respectively. Please note that the code, used by \cite{DS21} to obtain the true values by simulations, contained an error, leading to wrong values.

The solid black curve corresponds to the classical empirical estimator without PCA, the blue curves to the estimators based on $\hat H_{n,k}^{PCA}$ defined in \eqref{eq:hatHnkPCA}, and the red curves to the estimators using $\hat H_{n,k,\tilde k}^{PCA}$ defined in \eqref{eq:hatHnkktildePCA}, with $\tilde k=10$. Thereby, the colored dashed lines indicate the performance of the estimators based on PCA projections onto $p$ dimensional subspaces, while the solid lines refer to the case when the dimension is selected based on the $k$ (blue) respectively $\tilde k$ (red) observations with largest norms.

\begin{figure}
\includegraphics[width=0.9\textwidth]{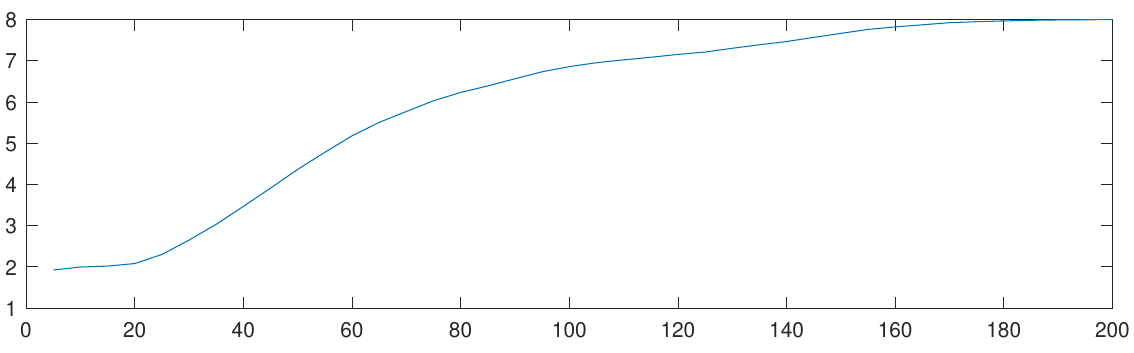}
\caption{Empirical mean of the selected dimension $\hat p$ versus $k$ in the Dirichlet model with $d=10$.}
 \label{fig:dimdistrdiri3_10_2}
\end{figure}

The estimators for the parameters (i) and (ii) based on $\hat H_{n,k,\tilde k}^{PCA}$ clearly perform best. Indeed, there is very little difference between the estimators with fixed dimension $p=2$ and with data-driven dimension $\hat p$, because in more than $99\%$ of the simulations the latter equals 2. The RMSE of these PCA-based estimators for the parameters  (iii) and (iv) is almost a constant function of $k$, which makes them easier to apply than the direct estimator without PCA, but for some values of $k$ the latter has a slightly smaller RMSE. The estimators using the variant of PCA that uses the $k$ most extreme observations (instead of merely $\tilde k=10$) usually does not perform well if $k$ is not very small. However, if the dimension is selected using the observations, then the performance of the estimators for (i)--(iii) is still slightly better than that of the direct estimator. Interestingly, this variant of the PCA-type estimators of (iii) and (iv) perform better than the one using the ``true'' dimension $p=2$. This can be explained by the fact that for larger values of $k$ the corresponding observations do not concentrate closely around a two-dimensional subspace. The selected dimension $\hat p$ is then substantially larger than 2, resulting in a smaller loss of information by the PCA projection. Figure \ref{fig:dimdistrdiri3_10_2} displays the average chosen dimension as a function of $k$; the distribution of $k$ is typically closely concentrated around this mean value, with more than 95\% of its mass concentrated in 2 or 3 values. Since for large values of $k$, the dimension reduction is not very pronounced, the performance of the estimators based on $\hat H_{n,k}^{PCA}$ is similar to the ones that do not use PCA at all. It then depends on the parameter to be estimated whether this is better or worse than first applying a PCA projection on a two-dimensional subspace which will often be quite different from the optimal subspace in the limit.

\begin{figure}
\includegraphics[width=0.9\textwidth]{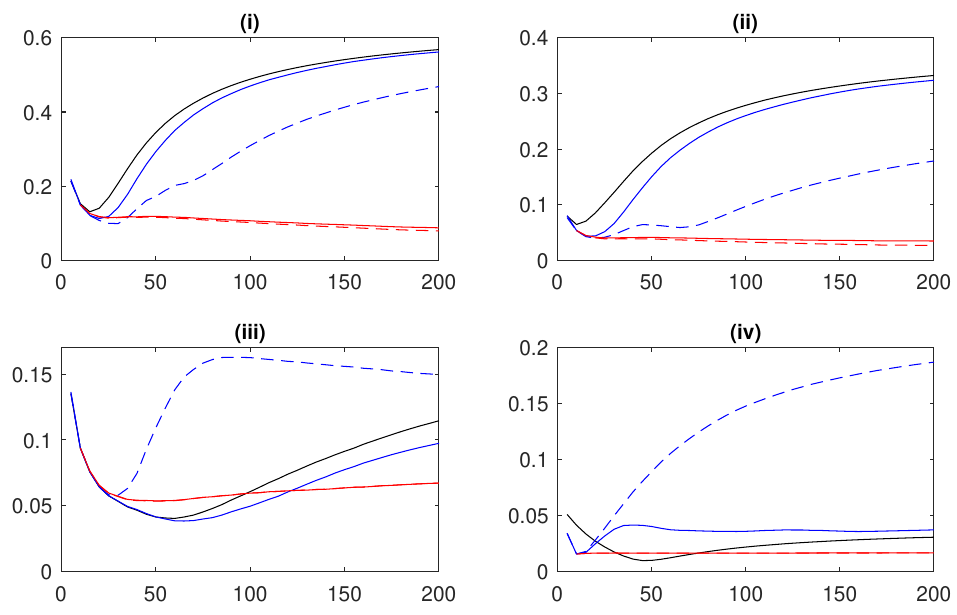}
\caption{RMSE of the estimators of the probabilities (i)--(iv) based on $\hat H_{n,k}$ (black, solid), $\hat H_{n,k}^{PCA}$ (blue, dashed) and $\hat H_{n,k,10}^{PCA}$ (red, dashed) versus $k$ for randomly rotated Dirichlet observations with fixed dimensions $p=2$ and $d=10$. The colored solid lines indicate the corresponding RMSE when the dimension is chosen as in \eqref{eq:choice_p}.}
 \label{fig:RMSEdiri_rot3_10_2}
\end{figure}

The performance of the estimators in the Dirichlet model with additional random rotation (resulting in the true values (i) 0.6527, (ii) 0.4016, (iii) 0.7618 and (iv) 0), shown in Figure \ref{fig:RMSEdiri_rot3_10_2}, is very similar. Despite the fact that the limit angular measure is no longer concentrated on a two-dimensional subspace, $\hat p$ equals 2 in more than $99\%$ of the simulations. This is probably due to the minor influence of the rotation on the reconstruction error and the choice of the parameter $\tau$ in the definition of $\hat p$, which tolerates a reconstruction error of 5\%.

\begin{figure}
\includegraphics[width=0.9\textwidth]{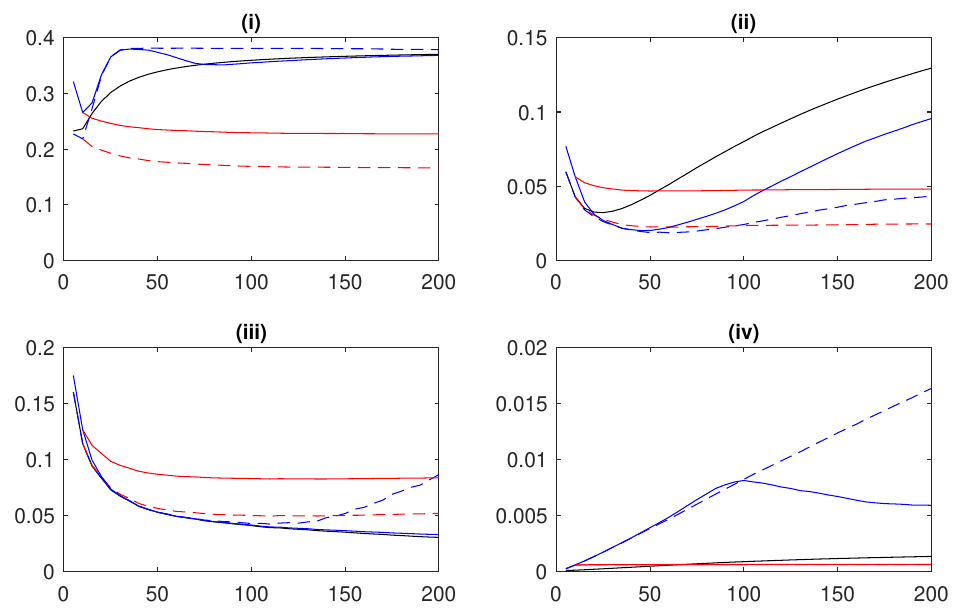}
\caption{RMSE of the estimators of the probabilities (i)--(iv) based on $\hat H_{n,k}$ (black, solid), $\hat H_{n,k}^{PCA}$ (blue, dashed) and $\hat H_{n,k,10}^{PCA}$ (red, dashed) versus $k$ in the Gumbel model with fixed dimensions $p=2$ and $d=10$. The colored solid lines indicate the corresponding RMSE when the dimension is chosen as in \eqref{eq:choice_p}.}
 \label{fig:RMSEGumbel_2_2_10_2}
\end{figure}

In contrast, in the Gumbel model with $\vartheta=2$, index of regular variation 2 and $k=10$, in about 93\%  of all simulations the chosen dimension equals the true value $p=2$, while in about 7\% it is $\hat p=1$. This leads to a substantial loss of accuracy of the estimators for the parameters (i)--(iii) compared with the estimators for fixed dimension $p=2$ (with $t_{(i)}=0.7$ and true values 0.3794, 0.2923 and $2^{-1/2}$, respectively), cf.\ Figure \ref{fig:RMSEGumbel_2_2_10_2}. However, the alternative version of the PCA procedure based on \eqref{eq:hatHnkktildePCA} and data driven choice of the dimension still performs much better than the direct estimator for parameter (i) and it is at least competitive for the parameters (ii) and (iv).

\subsubsection*{High dimensional models}

We now discuss the simulation results for the same models as before, but with overall dimension $d=100$ and angular measure concentrated on (the vicinity of) a subspace of dimension $p=5$. As one might have expected, it turns out that in the alternative approach based on \eqref{eq:hatHnkktildePCA} one should choose $\tilde k$ somewhat larger, because more observations are needed to estimate a subspace of higher dimension accurately. Here we show the results for $\tilde k=15$.

\begin{figure}
\includegraphics[width=0.9\textwidth]{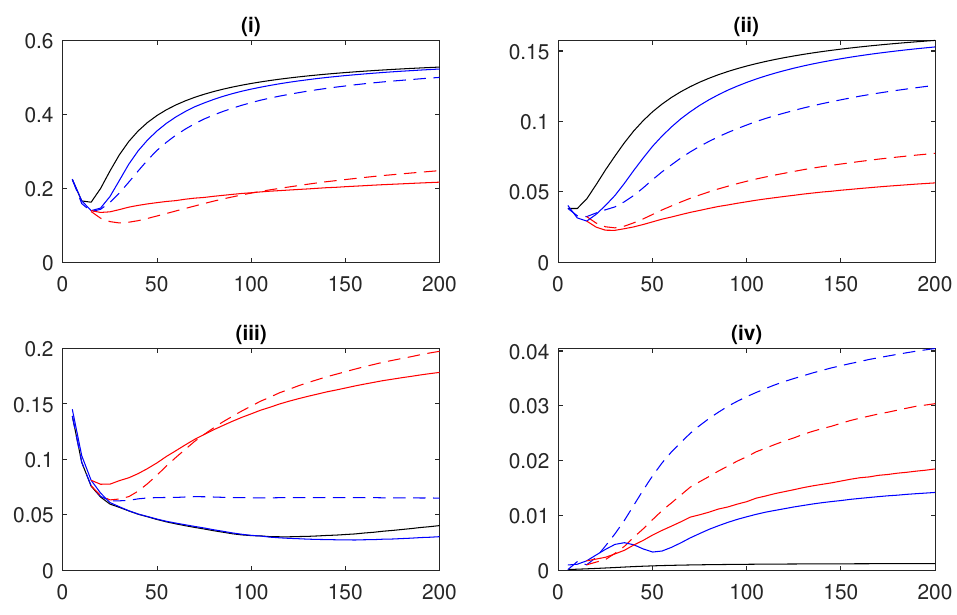}
\caption{RMSE of the estimators of the probabilities (i)--(iv) based on $\hat H_{n,k}$ (black, solid), $\hat H_{n,k}^{PCA}$ (blue, dashed) and $\hat H_{n,k,15}^{PCA}$ (red, dashed) versus $k$ in the Dirichlet model with fixed dimensions $p=5$ and $d=100$. The colored solid lines indicate the corresponding RMSE when the dimension is chosen as in \eqref{eq:choice_p}.}
 \label{fig:RMSEdiri3_100_5}
\end{figure}

Figure \ref{fig:RMSEdiri3_100_5} shows the RMSE of the estimator of the four parameters for the Dirichlet model (with $t_{(i)}=0.4$ and true values 0.5727, 0.1766, 0.5772 and 0). It is striking that in most cases the alternative PCA estimator with data driven choice of the dimension performs better than the one which uses the true dimension of the support of the angular measure, and for the parameters (iii) and (iv) this is also true for the PCA estimator based on \eqref{eq:hatHnkPCA}. Indeed, the chosen dimension based on the 15 vectors with largest norms equals $p=5$ in less than 20\% of the simulations, while in almost 78\% of the simulations $\hat p =4$. Somewhat surprisingly, for both statistics (ii) and (iv), for which the effect is strongest, projecting on a subspace of lower dimension reduces the bias, whereas the variance is reduced only for the estimator of (iv).

Compared with the direct estimator, the alternative PCA based estimators are clearly superior for (i) and (ii) and substantially worse for (iii) and (iv), but in the latter case all estimators are quite accurate. In contrast, the standard PCA based estimators projecting on a subspace of dimension $\hat p $  often performs better than the direct estimator, but the improvement is typically small.

The results for the Dirichlet model with additional rotation are quite similar and are thus not presented here. The main difference is that the PCA based procedures now performs better than the direct estimator of parameter (iv) if $k$ is less than 50.

For the Gumbel model (with $t_{(i)}=0.44$ and true parameters (i) 0.1835, (ii) 0.0827, (iii) $5^{-1/2}$ and (iv) 0), the alternative PCA based approach with data driven choice of $p$ performs slightly better than the one with fixed $p=5$ for the parameters (i) and (iv), but substantially worse for the probabilities (ii) and (iii). Here the chosen dimension is almost always equal to 3 or 4, with each of these values chosen with about the same probability. Since a projection on a three dimensional subspace cannot approximate the angular measure sufficiently well, as a result both the bias and the variance are substantially increased compared with the estimator with PCA projection on a five dimensional subspace.  In contrast, the standard PCA based estimators, which uses the same number of exceedances in each step of the procedures, now performs better if the dimension is chosen using the data for all parameters, but (ii). Moreover, the version with data driven choice of $p$ is never much worse than the direct estimator, and it is substantially better for probability (ii).

\begin{figure}
\includegraphics[width=0.9\textwidth]{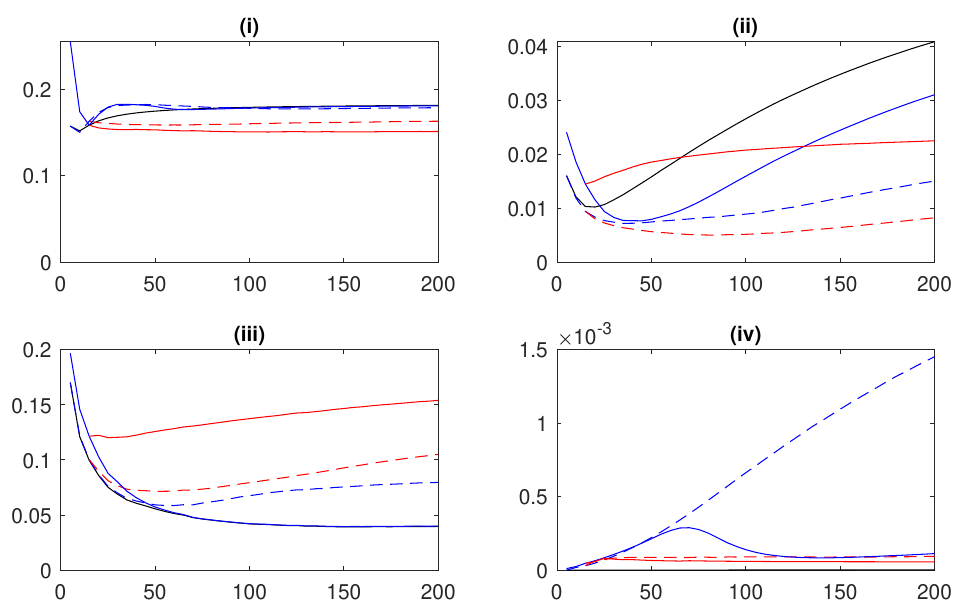}
\caption{RMSE of the estimators of the probabilities (i)--(iv) based on $\hat H_{n,k}$ (black, solid), $\hat H_{n,k}^{PCA}$ (blue, dashed) and $\hat H_{n,k,15}^{PCA}$ (red, dashed) versus $k$ in the Gumbel model with fixed dimensions $p=5$ and $d=100$. The colored solid lines indicate the corresponding RMSE when the dimension is chosen as in \eqref{eq:choice_p}.}
 \label{fig:RMSEGumbel_2_2_100_5}
\end{figure}

Overall, the picture not completely clear. However, it seems reasonable to use the alternative version of the PCA based estimators if the dimension of the data is moderate, and to perform PCA with the same number of exceedances for high-dimensional data. At least for the models considered here, this strategy almost never performs much worse than the direct estimator and in many cases it clearly outperforms the direct estimator.



\bibliographystyle{imsart-nameyear}
\bibliography{PCA_Lit}

\end{document}